\documentclass{amsart}

\usepackage{amsmath,amssymb,amsfonts,amsthm}
\usepackage{fullpage}
\usepackage{xurl}
\usepackage{hyperref}
\usepackage{xcolor}
\usepackage{enumitem}
\usepackage{tikz}
\usepackage{booktabs}
\usetikzlibrary{arrows.meta, positioning}

\newtheorem{theorem}{Theorem}
\newtheorem{lemma}[theorem]{Lemma}
\newtheorem{definition}[theorem]{Definition}

\newtheorem{proposition}[theorem]{Proposition}

\title{Erd\H{o}s Distance Results for $k$-Degenerate Graphs with $m$ Edges}

\author{Luv Barnwal}\email{l.barnwal@lse.ac.uk}
\address{Department of Mathematics, The London School of Economics and Political Science, London, United Kingdom}

\author{Yuda Chen}\email{math.yuda.chen@gmail.com}
\address{School of Mathematical Sciences and LPMC, Nankai University, Tianjin, China}

\author{Eon Lee}\email{eonlee1125@gmail.com}
\address{Department of Electrical Engineering and Computer Science, Gwangju Institute of Science and Technology, Gwangju, South Korea}

\author{Eyvindur Palsson}\email{palsson@vt.edu}
\address{Department of Mathematics, Virginia Tech, Blakcsburg VA, United States}

\author{Braiden Pulver}\email{bap028@bucknell.edu}
\address{Department of Mathematics and Statistics, Bucknell University, Lewisburg PA, United States}

\author{Steven Senger}\email{stevensenger@missouristate.edu}
\address{Department of Mathematics, Missouri State University, Springfield MO, United States}

\date{}

\begin{document}

\begin{abstract}
The Erd\H{o}s distinct distances problem asks how many distinct distances $N$ points must determine. In the plane, Guth and Katz determined the optimal power of $N$. A natural generalization fixes a finite graph $G$ with $m$ edges and records, over all injective embeddings of $G$ into a finite point set $P \subset \mathbb{R}^d$, the vectors of edge-lengths realized; we write $\Delta_G^\circ(P)$ for this set of edge-distance vectors. Taking $G$ to be a single edge recovers the classical problem. In the plane, the exponent $|V(G)|-1$ is known up to an arbitrarily small loss for every connected graph. A graph is $k$-degenerate if its vertices can be ordered so that each has at most $k$ neighbors among those preceding it. We show that if $G$ is $k$-degenerate, then for every $d \ge \max(2k,3)$ and every $\varepsilon > 0$, $|\Delta_G^\circ(P)| \gtrsim_{G,d,\varepsilon} N^{m/(d-1)-\varepsilon}$ for $|P|=N\ge |V(G)|$. For a connected graph with $m/(|V(G)|-1)>2(d-1)/d$, the exponent $m/(d-1)$ exceeds the spanning-tree exponent predicted by the pinned-distance conjecture. Our result is motivated by the threshold $\frac{d+k}{2}$ obtained by Borges, Foster, Ou, Palsson, and Romero Acosta for edge lengths of $k$-degenerate graphs in fractal sets. The proof combines a bisector-energy estimate for pinned $k$-stars with an entropy argument of Katz--Tardos type, applied along a construction order of $G$ and combined with an induction on the dimension.
\end{abstract}

\maketitle

\section{Introduction}
\label{sec:intro}

\subsection{History and background.}
Erd\H{o}s first asked for the minimum number of distinct distances $g(N)$ determined by $N$ points in the plane, and conjectured that $g(N)=\Omega(N/\sqrt{\log N})$~\cite{Erdos1946}. Erd\H{o}s himself proved $g(N)=\Omega(N^{1/2})$ by observing that the distances to two fixed points determine any other point up to two choices. In 2015, Guth and Katz~\cite{GuthKatz2015} proved $g(N)=\Omega(N/\log N)$ via polynomial partitioning, establishing the conjectured power of $N$. Their result built on earlier work: Moser~\cite{Moser} proved $g(N)=\Omega(N^{2/3})$; in 1992, Chung, Szemer{\'e}di and Trotter~\cite{Chung_Szemerédi_Trotter_1992} improved this to $\Omega(N^{4/5}/\log N)$; Sz\'ekely~\cite{Székely_1997} proved $\Omega(N^{4/5})$ in 1997; in 2001, Solymosi and T\'oth~\cite{Solymosi_Tóth_2001} proved $\Omega(N^{6/7})$; and in 2003, Tardos~\cite{Tardos_2003} improved this to $\Omega(N^{\frac{4e}{5e-1}-\varepsilon})$. In 2004, Katz and Tardos~\cite{KatzTardos2004} introduced an entropy inequality for the Erd\H{o}s distance problem, achieving a lower bound of $\Omega(N^{\frac{48-14e}{55-16e}}) \approx \Omega(N^{0.8641})$ distinct distances. 

In dimensions $d\ge4$, Solymosi and Vu~\cite{SolymosiVu} improved earlier estimates of Aronov, Pach, Sharir, and Tardos~\cite{ARONOV_PACH_SHARIR_TARDOS_2004}, obtaining the lower bound $\Omega_d(N^{2/d-2/(d(d+2))})$ for the total number of distinct distances. OpenAI~\cite{OpenAIHighDimDistances} subsequently proved the sharp bound $\Omega_d(N^{2/d})$ for every fixed $d\ge3$.

The \emph{pinned} problem asks how many different distances can be measured from a single point of the set. In the plane, OpenAI~\cite{OpenAIPlanarPinned} proved that for every $\delta>0$, all but $o(N)$ points of an $N$-point set determine at least $N^{1-\delta}$ distinct nonzero distances. In higher dimensions, the pinned estimates of Aronov et al.~\cite{ARONOV_PACH_SHARIR_TARDOS_2004} and the recurrences of Solymosi and Vu~\cite{SolymosiVu} give lower bounds for the number of distances visible from a suitable point.

Along with the distinct distance problem, Erd\H{o}s asked another important question, one which has had major breakthroughs recently. The question asks: "what is the maximum number of times the same distance can occur among $n$ points?". A suitably scaled grid gives at least $n^{1+c/\log\log n}$ unit distances, and Erd\H{o}s~\cite{Erdos1946} conjectured that this is essentially optimal, that is, at most $n^{1+O(1/\log\log n)}$. Spencer, Szemer{\'e}di, and Trotter~\cite{Spencer1984UnitDI} later produced a better upper bound, $O(n^{4/3})$. OpenAI~\cite{OpenAI} constructed, for infinitely many $n$, planar sets of $n$ points determining at least $n^{1+\delta}$ unit-distance pairs for a fixed $\delta>0$, disproving the conjecture. In subsequent work, OpenAI~\cite{OpenAIUnitUpper} improved the Spencer--Szemer\'edi--Trotter upper bound to $O(n^{4/3-\sigma})$ for an absolute constant $\sigma>0$. See further developments on quantifying the $\delta$ and $\sigma$ in \cite{ChenWang}.

One can replace a single edge by a finite graph $G$ and record the distances on its edges. Rudnev~\cite{Rudnev2012triangles} studied triangles in the plane. Iosevich and Passant~\cite{IosevichPassant2018} obtained sharp bounds for rigid frameworks in the plane. Dewar et al.~\cite[Corollary~1.8]{DFMNPW} proved that the measurement set has size $\gtrsim_G |P|^{|V(G)|-1}$ for every graph on at least three vertices that is generically rigid in the plane. For chains and graphs with a Hamiltonian path, Passant~\cite{Passant2020} obtained the lower bound $|P|^{|V(G)|-1}/\operatorname{polylog}|P|$ for the number of edge-length vectors. For every connected graph, the weak pinned planar distance theorem now gives the exponent $|V(G)|-1$ up to an arbitrarily small loss by the spanning-tree argument below.

Borges, Foster, Ou, Palsson, and Romero Acosta~\cite{BorgesFosterOuPalssonRomero2026} studied the corresponding fractal problem. For a $k$-degenerate graph with $m$ edges and $d>k$, they showed that if a compact set $E\subset\mathbb R^d$ satisfies $\dim_H(E)>(d+k)/2$, then the set of edge-length vectors realized in $E$ has positive $m$-dimensional Lebesgue measure. Their proof uses pinned $k$-stars. In the finite setting we count the different edge-length vectors of injective placements. Our lower bound depends on the number of edges and applies in every dimension $d\ge\max(2k,3)$.

\subsection{Main results and comparison with existing bounds.}

Throughout, $G=(V,E)$ denotes a finite graph with $n$ vertices and $m\ge1$ edges, and $P\subset\mathbb R^d$ a set of $N$ points. Fix an ordering of $E$. The set of edge-length vectors realized by injective placements of the vertices is
\[
\Delta_G^\circ(P)=\left\{\left(|\phi(u)-\phi(v)|\right)_{uv\in E}:\phi\colon V\hookrightarrow P\right\}.
\]
The coordinates follow the chosen edge order. Our main result is the following.

\begin{theorem}[see Theorem~\ref{thm:main}]\label{thm:intro-main}
Let $G=(V,E)$ be a finite graph with $m$ edges, and suppose $G$ is $k$-degenerate for some $k\ge1$. Then for every $d\ge\max(2k,3)$, every finite $P\subset\mathbb R^d$ with $|P|=N\ge|V|$, and every $\varepsilon>0$,
\[
|\Delta_G^\circ(P)|\gtrsim_{G,d,\varepsilon} N^{\frac{m}{d-1}-\varepsilon}.
\]
\end{theorem}

Another bound, valid in every dimension, follows from the spanning-tree argument of Iosevich and Passant~\cite{IosevichPassant2018}. They showed that in the plane the Erd\H{o}s pinned-distance conjecture implies lower bounds for all connected graphs. Their argument can be extended to $\mathbb R^d$ and to any pinned-distance estimate.

\begin{proposition}\label{prop:spanning}
Suppose there exist $\alpha,c>0$ such that every sufficiently large finite $Q\subset\mathbb R^d$ contains a point $x$ with $|\{|x-y|:y\in Q\}|\ge c|Q|^{\alpha}$. Then for every connected graph $G$ on $n$ vertices and every $P\subset\mathbb R^d$ with $|P|=N\ge n$,
\[
|\Delta_G^\circ(P)|\gtrsim N^{\alpha(n-1)},
\]
where the implicit constant may depend on $G$, $d$, and the pinned-distance estimate in the hypothesis.
\end{proposition}

For fixed $d$, the same $c$ and $\alpha$ work for every sufficiently large $Q$, although the point $x$ may depend on $Q$. The exponent $\alpha$ measures how the number of distances visible from one point grows with $|Q|$. This uniformity lets us apply the estimate to successive subsets of $P$ when embedding a spanning tree.

\begin{proof}
Root a spanning tree of $G$ and let $t$ be its depth. Given a sufficiently large set $Q$ of size $M$, at least half of its points see $c(M/2)^\alpha$ or more distinct distances in $Q$. Indeed, if more than half failed this bound, the hypothesis applied to the set of those points would give a contradiction. Starting with $P_t=P$, take $P_i$ to be the set of such points in $P_{i+1}$ for $i=t-1,\ldots,0$. Then $|P_i|\ge2^{-t}N$, and every point in $P_i$ sees $\gtrsim N^\alpha$ distances in $P_{i+1}$.

Place the root in $P_0$ and place the other vertices in order of increasing depth. If the parent lies in $P_i$, excluding all previously used points from $P_{i+1}$ removes at most $n$ possible lengths. Thus each partial placement extends in $\gtrsim N^\alpha$ different lengths on the next tree edge for all sufficiently large $N$. Keep one placement for each length at each step. The first edge at which two choices differ distinguishes the resulting vectors, giving $\gtrsim N^{\alpha(n-1)}$ vectors. For the bounded range of smaller $N\ge n$, at least one injective placement exists.
\end{proof}

In the plane, apply Proposition~\ref{prop:spanning} with $\alpha=1-\delta$ from the weak pinned planar distance theorem~\cite{OpenAIPlanarPinned}. For every fixed connected $G$ on $n\ge2$ vertices and every $\varepsilon>0$, choosing $\delta$ sufficiently small gives $|\Delta_G^\circ(P)|\gtrsim_{G,\varepsilon}N^{n-1-\varepsilon}$. If $G$ has $b$ connected components, partition $P$ into $b$ disjoint sets of comparable size and embed one component in each. The component vectors can be chosen independently, giving the exponent $n-b-\varepsilon$ after distributing the loss among the components. To see sharpness of the power, take $N$ equally spaced points on a line. In each component, the positions relative to one vertex determine its edge lengths, giving at most $O_G(N^{n-b})$ vectors in total.

For $d\ge3$, the pinned recurrences of Solymosi and Vu~\cite[Corollaries~2.3, 2.6 and Fact~2.8]{SolymosiVu}, started with the planar exponent $1-\delta$, give every exponent below
\[
\beta_d=
\begin{cases}
\dfrac{2(d+1)}{d^2+2d-2}, & d\ \text{even},\\[4pt]
\dfrac{2d^2}{d^3+d^2-3d+3}, & d\ \text{odd}.
\end{cases}
\]
For even $d$, this follows from their formula~(20) with initial dimension $2$ as $\delta\downarrow0$. For odd $d$, apply their one-step recurrence $\alpha\mapsto d\alpha/(d\alpha+d-1)$ from the preceding even dimension, taking $\beta_2=1$. In particular, $\beta_3=3/5$ and $\beta_4=5/11$. Proposition~\ref{prop:spanning} therefore gives every exponent below $(n-1)\beta_d$ for a connected graph.

The first edge of a spanning tree can also be chosen using the global distance bound~\cite{OpenAIHighDimDistances}. Root the tree at an edge, assigning depth zero to both endpoints, and form the sets $P_i$ as in the proof of Proposition~\ref{prop:spanning}. Since $|P_0|\gtrsim_G N$, pairs in $P_0$ realize $\gtrsim_d |P_0|^{2/d}$ different lengths. Keep one ordered pair for each length and extend it along the remaining $n-2$ tree edges. Each extension uses a pinned exponent arbitrarily close to $\beta_d$. Hence $|\Delta_G^\circ(P)|\gtrsim_{G,d,\varepsilon}N^{2/d+(n-2)\beta_d-\varepsilon}$. The conjectured pinned exponent $2/d$ in dimension $d$ would give the spanning-tree exponent $2(n-1)/d$. For trees this last power is sharp on the integer grid, which has only $O_d(N^{2/d})$ choices for each edge length.

For a connected graph with $n\ge2$ in the dimension range of Theorem~\ref{thm:intro-main}, write
\[
 A=\frac{m}{d-1},\qquad S=(n-1)\beta_d,\qquad
 U=\frac2d+(n-2)\beta_d,\qquad C=\frac{2(n-1)}d .
\]
These are exponents, with arbitrarily small losses understood. We have $S\le U\le C$. For trees, $\beta_d>1/(d-1)$, so even $S$ exceeds $A$. In particular, when $d=3$ the main theorem applies only to forests, for which the tree estimate applies componentwise. For graphs satisfying
\[
 \frac{m}{n-1}>\frac{2(d-1)}d ,
\]
the exponent $A$ exceeds even $C$. This includes planar triangulations and complete graphs on at least four vertices whenever the dimension hypothesis holds.

The effect of a single additional edge is visible already in $\mathbb R^4$. For the four-cycle $C_4$, the main exponent is $4/3$, while the two unconditional tree estimates give $S=15/11$ and $U=31/22$. Adding a diagonal gives the $2$-degenerate graph $K_4-e$ and raises the main exponent to $5/3$. The tree estimates have the same exponents as before, and even $C=3/2$ is smaller than $5/3$. Table~\ref{tab:comparison} gives the comparison for several graph families.

\begin{table}[ht]
\centering
\small
\setlength{\tabcolsep}{4pt}
\renewcommand{\arraystretch}{1.2}
\begin{tabular}{@{}lccccc@{}}
\toprule
Graph & Dimension & $A$ & $S>A$ & $U>A$ & $C>A$ \\
\midrule
Tree & $d\ge3$ & $(n-1)/(d-1)$ & always & always & always \\
Cycle $C_n$ & $d\ge4$ & $n/(d-1)$
  & \shortstack{$n\ge4$\\or $d\ge6$} & \shortstack{except\\$(3,4)$} & \shortstack{except\\$(3,4)$} \\
Planar triangulation & $d\ge10$ & $(3n-6)/(d-1)$ & never & never & never \\
$K_n$ & $d\ge 2n-2$ & $n(n-1)/(2(d-1))$ & never & never & never \\
\bottomrule
\end{tabular}
\caption{Comparison with the main exponent $A=m/(d-1)$. The exponents $S$ and $U$ are unconditional, and $C$ assumes the pinned-distance conjecture. For cycles $n\ge3$, the exception denotes $(n,d)=(3,4)$, where $S<U<A=C$. Planar triangulations and complete graphs have $n\ge4$.}
\label{tab:comparison}
\end{table}

\subsection{Proof methods.}
\label{subsec:methods}

Two methodological threads, both originating in the study of the classical
Erd\H{o}s distance problem, are combined and extended here. The argument uses two tools: one to control a single star, and one to combine stars across the graph.

The role of repeated perpendicular bisectors is implicit in Sz\'ekely's crossing-number argument~\cite{Székely_1997}. Solymosi and T\'oth~\cite{Solymosi_Tóth_2001} explicitly counted pairs with a common bisector, and Katz and Tardos~\cite{KatzTardos2004} developed this approach further. Lund, Sheffer, and de Zeeuw~\cite{LundShefferDeZeeuw2016} later defined and bounded the \emph{bisector energy}, the number of quadruples $(a,b,c,d)\in P^4$ with $a\ne b$ and $c\ne d$ for which $ab$ and $cd$ have the same perpendicular bisector. Lund~\cite{LundRefined2020}, Lund and Petridis~\cite{LundPetridis2020}, and Murphy, Rudnev, and Stevens~\cite{MurphyRudnevStevens2019} refined and extended their method.
Our popular star estimate (Theorem~\ref{thm:popstar}) is the $k$-anchor
generalization: we bound the energy $E(X)=\sum_t r_X(t)^2$ of a pinned $k$-star
by counting how many points can share a common perpendicular bisector \emph{hyperplane}. This count uses an incidence bound of Lund~\cite{Lund2021}, a
higher-dimensional generalization of the Szemer\'{e}di--Trotter theorem~\cite{SzemerediTrotter1983}
via the nondegeneracy framework of Elekes--T\'oth~\cite{ElekesToth2005} on nondegenerate hyperplanes, those whose points are spread across them rather than concentrated on a lower-dimensional flat.

\textbf{The entropy method.} Katz and Tardos~\cite{KatzTardos2004} introduced an entropy
inequality for the Erd\H{o}s distance problem. If $x$ is chosen uniformly from $P$ and $Z=|x-a|$, then $Z$ takes at most $D_a$ values, where $D_a$ is the number of distinct distances from $a$ to points of $P$. Hence $\log D_a\ge H(Z)$, and Jensen's inequality gives $H(Z)\ge\log(N^2/E(a))$ for the pinned energy $E(a)$. We generalize this from a single anchor to an arbitrary $k$-tuple (Lemma~\ref{lem:entropy_bound}), and then apply it along a fixed ordering $\mathcal O$ of the vertices of $G$. Placing the vertices onto distinct points of $P$ one at a time, the entropy contributed at each step is that of the new vertex’s distances to its already-placed neighbors; this is a pinned star, and the popular star estimate bounds it. Summing these contributions gives
$\log|\Delta_G^\circ(P)|\ge(\Phi_d(\mathcal O)-\varepsilon)\log N$ for an entropy
profile $\Phi_d$ defined from the order (Theorem~\ref{thm:entropy});
in dimensions $d=1,2$ no case split is needed.

The new ingredient is an induction on the dimension: when the point set $P$ has more than $N^\theta$ points on some $(d-2)$-flat, so that the popular star estimate does not apply, we restrict to the points on that flat and apply the entropy graph theorem in dimension at most $d-2$.
This is similar to the recursion of Solymosi and Vu~\cite[Theorem~2.2]{SolymosiVu}, which likewise reduces a $d$-dimensional distinct-distance problem to dimension $d-2$; our contribution is
to run this recursion through the entropy framework rather than on distinct
distances directly. The entropy profile $\Phi_r(\mathcal O)$ of Section~\ref{sec:entropy-profile} records the resulting exponent. It is non-increasing in $r$ (Lemma~\ref{lem:mono}), so a flat of dimension less than $d-2$ gives at least the bound for dimension $d-2$.

\subsection{Organization of the paper.}

In Section~\ref{sec:defs} we introduce the notation used throughout the paper. Given an ordered tuple $X=(a_1,\dots,a_k)$ of distinct points of $P$, we consider the set $S_X(P)$ of all vectors $(|x-a_1|,\dots,|x-a_k|)$ with $x\in P$, and call $a_1,\dots,a_k$ the \emph{anchors}. The \emph{energy} $E(X)$ is the number of pairs $(x,y)\in P^2$ whose vectors of distances to the anchors coincide. By the Cauchy--Schwarz inequality, $|S_X(P)|\ge N^2/E(X)$, so upper bounds on energy give lower bounds on $|S_X(P)|$. We write $\lambda(P)$ for the largest number of points of $P$ contained in a single affine $(d-2)$-flat. We write $k_i$ for the number of neighbors already placed when the $i$th vertex is placed.

In Section~\ref{sec:incidence} we state a theorem of Lund~\cite{Lund2021} that bounds the number of hyperplanes containing at least $s$ points of a finite set $A\subset\mathbb R^d$. The count includes only hyperplanes $\pi$ in which no affine $(d-2)$-flat contains more than half of the points of $A\cap\pi$. Using this bound, we prove upper bounds on $E(X)$ for $P\subset\mathbb R$ and $P\subset\mathbb R^2$.

In Section~\ref{sec:popstar} we prove the following estimate for $d\ge3$, which we call the popular star estimate (Theorem~\ref{thm:popstar}). Suppose that $\lambda(P)\le N^\theta$ for some $0\le\theta<1$. Then for every fixed $\eta>0$, all but $\eta N^k$ of the $k$-tuples $X$ of distinct points of $P$ satisfy
\[E(X)\lesssim N^{2-e_{d,k}(\theta)+\varepsilon},\qquad
e_{d,k}(\theta)=\begin{cases}\min\{1,k(1-\theta),k/(d-1)\},&1\le k<d,\\
1,&k\ge d .\end{cases}\]

In Section~\ref{sec:entropy-profile} we define, for each construction order $\mathcal O$ and each $d$, a quantity $\Phi_d(\mathcal O)$, which is the exponent that our argument produces in dimension $d$. For $d\ge3$ it is defined recursively in terms of $\Phi_{d-2}(\mathcal O)$. For $\theta\in[0,1)$, consider the following two exponents:
\begin{enumerate}
    \item the exponent $\sum_i e_{d,k_i}(\theta)$, which the estimate of Section~\ref{sec:popstar} yields when $\lambda(P)\le N^\theta$;
    \item the exponent $\theta\,\Phi_{d-2}(\mathcal O)$, which the bound in dimension $d-2$ yields when it is applied to the more than $N^\theta$ points of $P$ lying on a single $(d-2)$-flat.
\end{enumerate}
Then $\Phi_d(\mathcal O)$ is the supremum over $\theta$ of the smaller of these two exponents. We show that $\Phi_d(\mathcal O)$ is non-increasing in $d$, and that, when $d>\max_i k_i$, the value $1/\Phi_d(\mathcal O)$ is an explicit function of $1/\Phi_{d-2}(\mathcal O)$.

In Section~\ref{sec:entropy-thm} we prove that $|\Delta_G^\circ(P)|\gtrsim N^{\Phi_d(\mathcal O)-\varepsilon}$ for every construction order $\mathcal O$ of $G$. We choose an injective map $\phi\colon V(G)\to P$ uniformly at random and consider the vector $(|\phi(u)-\phi(v)|)_{uv\in E(G)}$, whose Shannon entropy is at most $\log|\Delta_G^\circ(P)|$. We bound this entropy from below by applying the chain rule along $\mathcal O$, and argue by induction on $d$.

In Section~\ref{sec:main} we iterate the relation between $\Phi_d$ and $\Phi_{d-2}$ established in Section~\ref{sec:entropy-profile}, separately in each parity class of $d$. We show that $\Phi_d(\mathcal O)=m/(d-1)$ for all $d\ge\max(2k,3)$ whenever $k_i\le k$ for all $i$, which yields Theorem~\ref{thm:intro-main}.

\section{Definitions}\label{sec:defs}

Let $P\subset\mathbb R^d$ be finite with $|P|=N$. For an ordered $k$-tuple of distinct points $X=(a_1,\ldots,a_k)\in P^k$, define the \emph{pinned $k$-star distance set}
\[
S_X(P)=\left\{\left(|x-a_1|,\ldots,|x-a_k|\right):x\in P\right\},
\]
and its \emph{energy}
\[
E(X)=\left|\left\{(x,y)\in P^2: |x-a_i|=|y-a_i|\ \text{for every } i=1,\ldots,k\right\}\right|.
\]
Writing $r_X(t)=|\{x\in P: (|x-a_1|,\ldots,|x-a_k|)=t\}|$, we have $E(X)=\sum_t r_X(t)^2$, so by Cauchy--Schwarz
\begin{equation}\label{eq:CS}
|S_X(P)|\ge \frac{N^2}{E(X)}.
\end{equation}
For $A\subseteq P$ write $A^{\underline k}=\{(a_1,\ldots,a_k)\in A^k: a_i \text{ pairwise distinct}\}$.

For $d\ge 3$, the \emph{codimension-two concentration} of $A\subset\mathbb R^d$ is
\[
\lambda(A)=\max_{\Gamma}|A\cap\Gamma|,
\]
the maximum taken over affine $(d-2)$-flats $\Gamma$.

\begin{definition}[Construction order]
A \emph{construction order} $\mathcal O$ for $G=(V,E)$ is an ordering $V=\{v_1,\ldots,v_n\}$ of its vertices. For each $i$, set
\[
J_i=\{u\in\{v_1,\ldots,v_{i-1}\}: v_iu\in E\},\qquad k_i=|J_i|,
\]
the set of neighbors of $v_i$ that precede it in the ordering. Each edge contributes once to the sum of the $k_i$, at its later endpoint, so
\[
m=|E(G)|=\sum_{i=1}^n k_i .
\]
\end{definition}

\begin{definition}[Degeneracy]
$G$ is \emph{$k$-degenerate} if it admits a construction order with $k_i\le k$ for every $i$. We call such an order a \emph{$k$-degenerate construction order}.
\end{definition}

Throughout, we assume that $G$ has at least one edge, so $m\ge1$. In particular, every $k$-degenerate construction order of $G$ has $k\ge1$. Additionally, we assume that $N\ge n=|V(G)|$; otherwise there is no injective map $V\to P$ and $\Delta_G^\circ(P)=\emptyset$.

\section{An incidence bound and the low-dimensional base cases}\label{sec:incidence}

An equality $|x-a|=|y-a|$ with $x\ne y$ puts the anchor $a$ on the perpendicular bisector of $xy$. To count these equalities, we use the following bound on the number of hyperplanes containing many points.

\begin{quote}
\textbf{Incidence bound} (Lund~\cite{Lund2021})\textbf{.}
Let $A\subset\mathbb R^d$, $|A|=N$. A hyperplane $\pi$ is \emph{$1/2$-nondegenerate} with respect to $A$ if no affine $(d-2)$-flat contains more than half the points of $A\cap\pi$. The number of $1/2$-nondegenerate hyperplanes containing at least $s\ge1$ points of $A$ is
\[
O_d\left(\frac{N^d}{s^{d+1}}+\frac{N^{d-1}}{s^{d-1}}\right).
\]
\end{quote}

\begin{lemma}[Linear base case]
\label{lem:d1}
Let $P\subset\mathbb R$, $|P|=N$, $k\ge1$. For every $X\in P^{\underline k}$, $E(X)\le 2N$, hence $|S_X(P)|\ge N/2$.
\end{lemma}
\begin{proof}
Fix $X=(a_1,\dots,a_k)$ and $t=(t_1,\dots,t_k)$. If $r_X(t)\neq 0$, then in particular any $x$ counted by $r_X(t)$ must satisfy
\[
|x-a_1|=t_1.
\]
In $\mathbb R$ this equation has at most two solutions, namely $x=a_1\pm t_1$. Therefore $r_X(t)\le 2$ for every $t$. Since $\sum_t r_X(t)=N$,
\[
E(X)=\sum_t r_X(t)^2\le 2\sum_t r_X(t)=2N.
\]
Now apply \eqref{eq:CS}.
\end{proof}

\begin{lemma}[Planar base case]
\label{lem:d2}
Let $P\subset\mathbb R^2$, $|P|=N$.
\begin{enumerate}
\item[(a)] For $k\ge2$ and every $X\in P^{\underline k}$, $E(X)\le 2N$.
\item[(b)] For $k=1$, every $\eta,\varepsilon>0$: all but $\eta N$ points $a\in P$ satisfy $E(a)\lesssim_{\eta,\varepsilon}N^{4/3+\varepsilon}$, hence $|S_a(P)|\gtrsim_{\eta,\varepsilon}N^{2/3-\varepsilon}$.
\end{enumerate}
\end{lemma}
\begin{proof}
\medskip\noindent\textbf{Part (a).}
Let $X=(a_1,\dots,a_k)$ with $k\ge2$. For any $t$, any $x$ counted by $r_X(t)$ must satisfy
\[
|x-a_1|=t_1,\qquad |x-a_2|=t_2.
\]
Since $X\in P^{\underline k}$ has distinct entries, $a_1\ne a_2$, so the two circles have distinct centers and are therefore distinct; two distinct circles intersect in at most two points. Hence $r_X(t)\le2$ for every $t$, so $E(X)\le2N$.

\medskip\noindent\textbf{Part (b).}
Summing over anchors,
\[
\sum_{a\in P} E(a) \le N\cdot N + \sum_\ell \mu(\ell)|P\cap\ell|,
\]
where $\ell$ ranges over perpendicular bisector lines and $\mu(\ell)$ is the multiplicity of $\ell$. Dyadically, let $S_m=\{\ell: m\le |P\cap\ell|<2m\}$. For every $m$, the trivial bound $\sum_\ell\mu(\ell)\le N^2$ gives $\sum_{\ell\in S_m}\mu(\ell)|P\cap\ell|\lesssim N^2m$. For $m\ge2$, the Szemer\'edi--Trotter theorem~\cite{SzemerediTrotter1983} gives
\[
\sum_{\ell\in S_m}\mu(\ell)|P\cap\ell|
\lesssim N\cdot m\left(\frac{N^2}{m^3}+\frac Nm\right)
= \frac{N^3}{m^2}+N^2.
\]
Hence each class $S_m$ contributes at most $\min(N^2m,\,N^3/m^2)+N^2$. Since
\[
\min\Bigl(N^2m,\ \frac{N^3}{m^2}\Bigr)\le (N^2m)^{2/3}\Bigl(\frac{N^3}{m^2}\Bigr)^{1/3}=N^{7/3},
\]
each class contributes $O(N^{7/3})$. Summing over the $O(\log N)$ dyadic classes,
\[
\sum_{a\in P}E(a)\lesssim N^{7/3+o(1)},
\qquad \mathbb E_a[E(a)]\lesssim N^{4/3+o(1)}.
\]
By Markov's inequality, all but an $\eta$-fraction of $a\in P$ satisfy $E(a)\lesssim_{\eta,\varepsilon}N^{4/3+\varepsilon}$; apply \eqref{eq:CS}.
\end{proof}

\section{The popular star estimate in dimension \texorpdfstring{$d\geq3$}{d≥3}}\label{sec:popstar}

For $d\ge3$, $1\le k\le d-1$, $0\le\theta<1$, define
\[
e_{d,k}(\theta)=\min\left(1,\ k(1-\theta),\ \frac{k}{d-1}\right),
\]
and $e_{d,k}(\theta)=1$ for $k\ge d$, $e_{d,0}(\theta)=0$.

\begin{theorem}[Popular star estimate]\label{thm:popstar}
Let $P\subset\mathbb R^d$, $|P|=N$, $d\ge3$, and suppose $\lambda(P)\le N^\theta$ for some $0\le\theta<1$. Let $k\ge1$, $\eta\in(0,1)$ and $\varepsilon>0$.
\begin{enumerate}
\item[(i)] If $k<d$, then all but at most $\eta N^k$ tuples $X\in P^{\underline k}$ satisfy
\[
|S_X(P)|\gtrsim_{\eta,d,k,\varepsilon} N^{e_{d,k}(\theta)-\varepsilon}.
\]
\item[(ii)] If $k\ge d$, then there is $N_0=N_0(d,k,\theta,\eta)$ such that for $N\ge N_0$, all but at most $\eta N^k$ tuples $X\in P^{\underline k}$ satisfy $E(X)\le 2N$, and hence $|S_X(P)|\ge N/2$.
\end{enumerate}
In both cases, all but at most $\eta N^k$ tuples $X\in P^{\underline k}$ satisfy
\[
E(X)\lesssim_{\eta,d,k,\varepsilon,\theta} N^{2-e_{d,k}(\theta)+\varepsilon}.
\]
For (ii) and $N<N_0$, this follows from the trivial bound $E(X)\le N^2\le N_0N$.
\end{theorem}
\begin{proof}
\medskip\noindent\textbf{Case 1 ($1\le k<d$).}
Summing $E(X)$ over ordered distinct $k$-tuples, the diagonal $x=y$ contributes at most $N^{k+1}$. For $x\ne y$, let $B(x,y)$ be the perpendicular bisector hyperplane of $xy$. Let
\[
\mu(\pi)=|\{(x,y)\in P^2: x\ne y,\ B(x,y)=\pi\}|.
\]
Then, the off-diagonal contribution to $\sum_X E(X)$ is bounded by $\sum_\pi \mu(\pi)|P\cap\pi|^k$.

We dyadically decompose this sum over sets $H_m = \{\pi : m \le |P\cap\pi| < 2m\}$. Let
\[
\Sigma_m = \sum_{\pi \in H_m} \mu(\pi)|P\cap\pi|^k \lesssim m^k \sum_{\pi \in H_m} \mu(\pi).
\]
We use two bounds for $\sum_{\pi\in H_m}\mu(\pi)$:
\begin{itemize}
    \item \emph{Trivial bound:} $\sum_{\pi \in H_m} \mu(\pi) \le N^2$.
    \item \emph{Incidence bound:} Since $\mu(\pi)\le N$, $\sum_{\pi \in H_m} \mu(\pi) \le N |H_m|$. For $m>2\lambda(P)$, any $\pi \in H_m$ is $1/2$-nondegenerate. By Lund's incidence bound, $|H_m| \lesssim_d N^d/m^{d+1} + N^{d-1}/m^{d-1}$. Thus $\sum_{\pi \in H_m} \mu(\pi) \lesssim N^{d+1}/m^{d+1} + N^d/m^{d-1}$.
\end{itemize}
The scales $m\le2\lambda(P)$ contribute $O_k(N^2\lambda(P)^k)$ in total, since $|P\cap\pi|<4\lambda(P)$ on these scales and $\sum_\pi\mu(\pi)\le N^2$. For $m>2\lambda(P)$, we have $\Sigma_m \lesssim m^k \min\left(N^2,\ \frac{N^{d+1}}{m^{d+1}}+\frac{N^d}{m^{d-1}}\right)$.

Since $\min(A,B+C)\le\min(A,B)+\min(A,C)$,
\[
\Sigma_m \lesssim m^k \min\left(N^2,\ \frac{N^{d+1}}{m^{d+1}}\right) + m^k \min\left(N^2,\ \frac{N^d}{m^{d-1}}\right).
\]
We bound the two terms separately.
\begin{itemize}
    \item \emph{Term 1:} The value of $m$ where the two terms are balanced ($N^2 = \frac{N^{d+1}}{m^{d+1}}$) is $m_1 = N^{\frac{d-1}{d+1}}$. For scales $m \le m_1$, the trivial bound $N^2 m^k$ increases with $m$. For scales $m > m_1$, the incidence bound $\frac{N^{d+1}}{m^{d+1}} m^k = N^{d+1} m^{k-d-1}$ decreases with $m$ because $k \le d-1 < d+1$. Therefore, the maximum occurs at the crossover $m = m_1$. Summing over dyadic scales yields $\lesssim N^2 m_1^k \log N = N^{2 + k\frac{d-1}{d+1}} \log N$.
    \item \emph{Term 2:} The value of $m$ where the two terms are balanced ($N^2 = \frac{N^d}{m^{d-1}}$) is $m_2 = N^{\frac{d-2}{d-1}}$. Summing over dyadic scales, the maximum contribution is similarly governed by $m_2$. If $k < d-1$, the sum is geometrically dominated by $m_2$. If $k = d-1$, all dyadic scales contribute equally, giving an extra factor $O(\log N)$. Both yield $\lesssim N^2 m_2^k \log N = N^{2 + k\frac{d-2}{d-1}} \log N$.
\end{itemize}
To determine the dominant term, observe $(d-1)^2 = d^2 - 2d + 1$ and $(d-2)(d+1) = d^2 - d - 2$. Subtracting yields $(d^2 - d - 2) - (d^2 - 2d + 1) = d - 3$. Thus, for all $d \ge 3$, $(d-1)^2 \le (d-2)(d+1)$, meaning $\frac{d-1}{d+1} \le \frac{d-2}{d-1}$.
Therefore Term 2 dominates, with equality when $d=3$. Since there are $O(\log N)$ dyadic scales and $\log N=N^{o(1)}$, we obtain
\[
\sum_\pi\mu(\pi)|P\cap\pi|^k\lesssim N^2\lambda(P)^k + N^{2+\frac{k(d-2)}{d-1}+o(1)}.
\]
If $N<k$, there are no distinct $k$-tuples. Otherwise, choose $X$ uniformly from $P^{\underline k}$ and write $e=e_{d,k}(\theta)$. Dividing the diagonal and off-diagonal bounds by $|P^{\underline k}|\gtrsim_k N^k$ and using $\lambda(P)\le N^\theta$, we obtain
\[
\mathbb E[E(X)]\le C_{d,k,\varepsilon}N^{2-e+\varepsilon/2},
\qquad 2-e=\max\{1,\,2-k+k\theta,\,2-k/(d-1)\}.
\]
Set $T=C_{d,k,\varepsilon}N^{2-e+\varepsilon/2}/\eta$. Markov's inequality for the nonnegative random variable $E(X)$ gives
\[
\mathbb P\bigl(E(X)>T\bigr)\le\frac{\mathbb E[E(X)]}{T}\le\eta.
\]
Thus at most $\eta|P^{\underline k}|\le\eta N^k$ tuples are exceptional. For every other tuple, \eqref{eq:CS} gives $|S_X(P)|\ge N^2/T\gtrsim_{\eta,d,k,\varepsilon}N^{e-\varepsilon/2}$, as required.

\medskip\noindent\textbf{Case 2 ($k\ge d$).}
If the first $d$ points of $X$ are affinely independent, then for any prescribed distance vector $t=(t_1,\dots,t_k)$ the equations
\[
|x-a_i|=t_i,\qquad i=1,\dots,d,
\]
after subtracting the first equation squared from each of the others, yield $d-1$ linear equations with normals $a_i-a_1$, $2\le i\le d$. These vectors are linearly independent, so the solution set is an affine line $L_t$. Intersecting $L_t$ with the first sphere ($\{x:|x-a_1|=t_1\}$) gives at most two points, so $r_X(t)\le2$ for every $t$. Therefore
\[
E(X)=\sum_t r_X(t)^2\le 2\sum_t r_X(t)=2N.
\]
It remains to control the tuples whose first $d$ points are not affinely independent. Let $2\le i\le d$ be the first index such that $a_1,\dots,a_i$ are affinely dependent. Then $a_i$ lies in the affine span of $a_1,\dots,a_{i-1}$, which has dimension at most $i-2$ and is therefore contained in some affine $(d-2)$-flat. Hence, once $a_1,\dots,a_{i-1}$ are fixed, there are at most $\lambda(P)\le N^\theta$ choices for $a_i$.
Summing over $i=2,\dots,d$ shows that the number of ordered $d$-tuples with a dependent prefix is at most $C_dN^{d-1+\theta}$, and extending the remaining $k-d$ coordinates freely gives at most $C_dN^{k-1+\theta}$ bad $k$-tuples. Since $\theta<1$, this is at most $\eta N^k$ when
\[
N\ge N_0:=\left(C_d/\eta\right)^{1/(1-\theta)}.
\]
For all remaining tuples, $E(X)\le2N$, which is the exponent $1=e_{d,k}(\theta)$. This completes the proof.
\end{proof}

\section{The entropy profile and its properties}\label{sec:entropy-profile}

The star estimate gives one exponent for each vertex with previously placed neighbors. We combine these exponents and compare their sum with the bound obtained by restricting to a $(d-2)$-flat.

\begin{definition}[Entropy profile]\label{def:phi}
For a construction order $\mathcal O$ of $G$ with degrees $(k_i)_{i=1}^n$, define
\[
\Phi_1(\mathcal O)=|\{i:k_i\ge1\}|,\qquad
\Phi_2(\mathcal O)=\frac23|\{i:k_i=1\}|+|\{i:k_i\ge2\}|.
\]
For $d\ge3$, define recursively
\[
E_d(\mathcal O,\theta)=\sum_{i=1}^n e_{d,k_i}(\theta),\qquad
\Phi_d(\mathcal O)=\sup_{0\le\theta<1}\min\bigl(E_d(\mathcal O,\theta),\ \theta\,\Phi_{d-2}(\mathcal O)\bigr).
\]
\end{definition}

\begin{lemma}[Monotonicity]\label{lem:mono}
$\Phi_r(\mathcal O)$ is non-increasing in $r$: for every $r\ge2$, $\Phi_r(\mathcal O)\le\Phi_{r-1}(\mathcal O)$.
\end{lemma}
\begin{proof}
We argue by strong induction on $r$.

\medskip\noindent\textbf{Base case ($r=2$).}
By definition, for each index $i$ with $k_i\ge1$, the term contributed to $\Phi_1$ is $1$, while the term contributed to $\Phi_2$ is $2/3$ if $k_i=1$ and $1$ if $k_i\ge2$; in either case it is $\le1$. Indices with $k_i=0$ contribute $0$ to both. Summing termwise over $i$ gives $\Phi_2(\mathcal O)\le\Phi_1(\mathcal O)$.

\medskip\noindent\textbf{Base case ($r=3$).}
Extend the definition of $e_{d,k}$ to $d=2$ by setting $e_{2,1}=2/3$, $e_{2,k}=1$ for $k\ge2$, and $e_{2,0}=0$, so that $\Phi_2(\mathcal O)=\sum_i e_{2,k_i}$ matches Definition~\ref{def:phi}. Then for every $\theta\in[0,1)$ and every $k\ge0$,
\[
e_{3,k}(\theta)\le e_{2,k}.
\]
Indeed, for $k=0$ both sides are $0$; for $k=1$,
\[
e_{3,1}(\theta)=\min(1,\,1-\theta,\,1/2)\le1/2<2/3=e_{2,1};
\]
for $k=2$,
\[
e_{3,2}(\theta)=\min(1,\,2(1-\theta))\le1=e_{2,2};
\]
and for $k\ge3$ both sides equal $1$. Summing over $i$ gives $E_3(\mathcal O,\theta)\le\Phi_2(\mathcal O)$ for every $\theta$, hence
\[
\Phi_3(\mathcal O)=\sup_{0\le\theta<1}\min\bigl(E_3(\mathcal O,\theta),\,\theta\Phi_1(\mathcal O)\bigr)\le\sup_{0\le\theta<1}E_3(\mathcal O,\theta)\le\Phi_2(\mathcal O).
\]

\medskip\noindent\textbf{Inductive step ($r\ge4$).}
Suppose, as the strong induction hypothesis, that $\Phi_s(\mathcal O)\le\Phi_{s-1}(\mathcal O)$ holds for every $s$ with $2\le s<r$. Since $r\ge4$, both $r-2$ and $r-3$ satisfy $2\le r-2<r$, so in particular the hypothesis applies at $s=r-2$:
\begin{equation}\label{eq:ih}
\Phi_{r-2}(\mathcal O)\le\Phi_{r-3}(\mathcal O).
\end{equation}

We claim that
\[
e_{r,k}(\theta)\le e_{r-1,k}(\theta)\qquad\text{for every }k\ge0,\ \theta\in[0,1).
\]
If $k=0$, both sides vanish. If $k\ge r-1$, then $e_{r-1,k}(\theta)=1$ by convention, while $e_{r,k}(\theta)\le1$ always, since it is either $1$ or a minimum one of whose arguments is $1$. Finally, if $1\le k\le r-2$, then $k<r-1<r$ and $r-1\ge3$, so both quantities are given by the min-formula:
\[
e_{r,k}(\theta)=\min\Bigl(1,\,k(1-\theta),\,\frac{k}{r-1}\Bigr),\qquad
e_{r-1,k}(\theta)=\min\Bigl(1,\,k(1-\theta),\,\frac{k}{r-2}\Bigr).
\]
Since $\frac{k}{r-1}\le\frac{k}{r-2}$ and the other two arguments coincide, and a minimum is non-decreasing in each argument, the claim follows.
Summing over $i$,
\begin{equation}\label{eq:Er}
E_r(\mathcal O,\theta)\le E_{r-1}(\mathcal O,\theta)\qquad\text{for every }\theta\in[0,1).
\end{equation}

Now fix any $\theta\in[0,1)$. Using \eqref{eq:Er} and monotonicity of $\min(\cdot,B)$ in its first argument,
\[
\min\bigl(E_r(\mathcal O,\theta),\,\theta\Phi_{r-2}(\mathcal O)\bigr)\le\min\bigl(E_{r-1}(\mathcal O,\theta),\,\theta\Phi_{r-2}(\mathcal O)\bigr).
\]
Using \eqref{eq:ih}, $\theta\ge0$, and monotonicity of $\min(A,\cdot)$ in its second argument,
\[
\min\bigl(E_{r-1}(\mathcal O,\theta),\,\theta\Phi_{r-2}(\mathcal O)\bigr)\le\min\bigl(E_{r-1}(\mathcal O,\theta),\,\theta\Phi_{r-3}(\mathcal O)\bigr).
\]
Chaining these,
\[
\min\bigl(E_r(\mathcal O,\theta),\,\theta\Phi_{r-2}(\mathcal O)\bigr)\le\min\bigl(E_{r-1}(\mathcal O,\theta),\,\theta\Phi_{r-3}(\mathcal O)\bigr)\qquad\text{for every }\theta\in[0,1).
\]
Taking the supremum over $\theta\in[0,1)$ on both sides preserves the inequality:
\[
\Phi_r(\mathcal O)=\sup_{0\le\theta<1}\min\bigl(E_r(\mathcal O,\theta),\,\theta\Phi_{r-2}(\mathcal O)\bigr)\le\sup_{0\le\theta<1}\min\bigl(E_{r-1}(\mathcal O,\theta),\,\theta\Phi_{r-3}(\mathcal O)\bigr)=\Phi_{r-1}(\mathcal O).
\]

This closes the induction, so $\Phi_r(\mathcal O)\le\Phi_{r-1}(\mathcal O)$ holds for every $r\ge2$.
\end{proof}

\begin{lemma}\label{lem:deg}
Let $\mathcal O$ be a $k$-degenerate construction order with $m=\sum_i k_i\ge1$. For every $r\ge\max(k+1,3)$, the profile satisfies
\[
\frac{1}{\Phi_r(\mathcal O)} = \max\left( \frac{r-1}{m},\ \frac{1}{\Phi_{r-2}(\mathcal O)} + \frac{1}{m} \right).
\]
\end{lemma}
\begin{proof}
We first note that $\Phi_s(\mathcal O)>0$ for every $s\ge1$. For $s=1,2$ this holds because $m\ge1$. For $s\ge3$, assume inductively that $\Phi_{s-2}>0$. Since some $k_i\ge1$, we have $E_s(\mathcal O,0)>0$, so by continuity $E_s(\mathcal O,\theta)>0$ for all small $\theta>0$. For such $\theta$, also $\theta\Phi_{s-2}>0$, hence $\Phi_s>0$.

Now let $r\ge\max(k+1,3)$. Since $r\ge3$, $\Phi_r$ is given by the recursive definition and $\Phi_{r-2}$ is defined. Since $r>k$, we have $k_i\le k\le r-1<r$ for all $i$, so every $e_{r,k_i}$ is given by the min-formula (or vanishes when $k_i=0$). As $k_i/(r-1)\le1$, it simplifies to
\[
e_{r,k_i}(\theta)=\min\Bigl(1,\,k_i(1-\theta),\,\frac{k_i}{r-1}\Bigr)=k_i\min\Bigl(1-\theta,\,\frac{1}{r-1}\Bigr).
\]
Summing over $i$ gives $E_r(\theta)=m\min\bigl(1-\theta,\frac{1}{r-1}\bigr)$, so
\[
\Phi_r=\sup_{0\le\theta<1}\min\Bigl(m(1-\theta),\ \frac{m}{r-1},\ \theta\Phi_{r-2}\Bigr).
\]
The function $\theta\mapsto\min\bigl(m(1-\theta),\frac{m}{r-1}\bigr)$ is constant at $\frac{m}{r-1}$ for $\theta\le\frac{r-2}{r-1}$ and then strictly decreases to $0$. Since $\Phi_{r-2}>0$, the function $\theta\mapsto\theta\Phi_{r-2}$ is strictly increasing. Their unique intersection determines the supremum.
\begin{itemize}
\item If they intersect on the plateau, then $\Phi_r=\frac{m}{r-1}$. This occurs exactly when $\frac{m}{r-1}\le\frac{r-2}{r-1}\Phi_{r-2}$.
\item Otherwise, they intersect on the slope, where $m(1-\theta)=\theta\Phi_{r-2}$. This gives $\Phi_r=\frac{m\Phi_{r-2}}{m+\Phi_{r-2}}$.
\end{itemize}
Taking reciprocals in both cases gives the stated maximum.
\end{proof}
\section{The entropy graph theorem}\label{sec:entropy-thm}

We place the vertices in construction order and estimate the entropy of the edge lengths added at each step. For a finite-valued random variable $Z$ with probabilities $p_z=\mathbb P(Z=z)$, its Shannon entropy is
\[
H(Z)=-\sum_{z:\,p_z>0}p_z\log p_z.
\]
The following lemma bounds this entropy from below when $Z$ records the distances to a fixed tuple of anchors.

\begin{lemma}
\label{lem:entropy_bound}
Let $X = (a_1, \dots, a_k) \in P^{\underline{k}}$ be a fixed tuple, and let $P' \subseteq P$ be a subset of points with $|P'| = M \ge 1$. If a point $x$ is chosen uniformly at random from $P'$, the Shannon entropy of the random distance vector $Z = (|x-a_1|, \ldots, |x-a_k|)$ satisfies
\[H(Z) \ge \log \frac{M^2}{E_P(X)},\]
where $E_P(X)$ is the energy of $X$ computed in $P$.
\end{lemma}
\begin{proof}
The random variable $Z$ takes value $z$ with probability $p_z = r'(z) / M$, where $r'(z)$ is the number of points in $P'$ yielding distance vector $z$. By Jensen’s inequality,
\[H(Z) = \sum_z p_z (-\log p_z) \ge -\log \Bigl( \sum_z p_z^2 \Bigr).\]
Substituting $p_z = r'(z)/M$, and using the fact that $r'(z) \le r_X(z)$, we obtain
\[\sum_z p_z^2 = \frac{1}{M^2} \sum_z r'(z)^2 \le \frac{1}{M^2} \sum_z r_X(z)^2 = \frac{E_P(X)}{M^2}.\]
Because $-\log$ is a monotonically decreasing function, $H(Z) \ge -\log(E_P(X)/M^2) = \log(M^2/E_P(X))$.
\end{proof}

\begin{theorem}\label{thm:entropy}
For every fixed graph $G$, construction order $\mathcal O$, dimension $d\ge1$, finite $P\subset\mathbb R^d$ with $|P|=N$, and $\varepsilon>0$,
\[|\Delta_G^\circ(P)|\gtrsim_{G,\mathcal O,d,\varepsilon} N^{\Phi_d(\mathcal O)-\varepsilon}.\]
\end{theorem}
\begin{proof}
We argue by strong induction on the ambient dimension $d$, with base cases $d=1,2$. Choose a map $\phi:V(G)\hookrightarrow P$ uniformly at random from all $\frac{N!}{(N-n)!}$ injective mappings, and set $Y_i=(|\phi(v_i)-\phi(u)|)_{u\in J_i}$, $Y=(Y_1,\ldots,Y_n)$.
Since $N\ge n$, some injective map exists, so $|\Delta_G^\circ(P)|\ge1$; hence it suffices to prove the bound for $N$ larger than a constant depending on $G$, $\mathcal O$, $d$ and $\varepsilon$.
The restriction of $\phi$ to any set of vertices is a uniform injective map, so the tuple $X_i=(\phi(u))_{u\in J_i}$ is uniform on $P^{\underline{k_i}}$.

The vector $Y$ lists each edge length once, at the later endpoint of the edge. Its coordinates are a fixed reordering of the edge-length vector, so it takes at most $|\Delta_G^\circ(P)|$ values. Thus
\[\log|\Delta_G^\circ(P)|\ge H(Y)=\sum_{i=1}^n H(Y_i\mid Y_1,\ldots,Y_{i-1}).\]
Let $\mathcal F_i=(\phi(v_1),\ldots,\phi(v_{i-1}))$. A \emph{history at step $i$} is a possible value $h=(p_1,\ldots,p_{i-1})$ of this tuple, so that $\phi(v_j)=p_j$ for $j<i$. The entries of $h$ are distinct, and at step $1$ the history is the empty tuple. Write $X_i(h)$ for the subtuple of $h$ indexed by the vertices in $J_i$, in construction order. Since $(Y_1,\ldots,Y_{i-1})$ is a function of $\mathcal F_i$,
\begin{equation}\label{eq:cond}
H(Y_i\mid Y_1,\ldots,Y_{i-1})\ \ge\ H(Y_i\mid \mathcal F_i)\ =\ \sum_{h}\mathbb P(\mathcal F_i=h)\,H(Y_i\mid \mathcal F_i=h).
\end{equation}
Conditionally on $\mathcal F_i=h$, the point $\phi(v_i)$ is uniform on the $N-i+1$ points of $P$ not used by $h$. Applying Lemma~\ref{lem:entropy_bound} with $P'$ equal to this set, which has $N-i+1\ge N-n+1\ge1$ points, gives
\begin{equation}\label{eq:pointwise}
H(Y_i\mid \mathcal F_i=h)\ \ge\ \log\frac{(N-i+1)^2}{E(X_i(h))}\ \ge\ \log\frac{(N-n+1)^2}{E(X_i(h))}.
\end{equation}
Steps with $k_i=0$ have $Y_i$ empty and do not contribute to $H(Y)$.

\medskip\noindent\textbf{Base case ($d=1$).}
The linear base case (Lemma~\ref{lem:d1}) gives $E(X_i)\le2N$ whenever $k_i\ge1$. Therefore, by \eqref{eq:pointwise}, $H(Y_i\mid\mathcal F_i=h)\ge\log(N/2)-O_G(1)$ for every history $h$, and by \eqref{eq:cond},
\[H(Y_i\mid Y_1, \ldots, Y_{i-1})\ge \log(N/2)-O_G(1).\]
Summing over the steps with $k_i\ge1$ gives $H(Y)\ge\Phi_1(\mathcal O)\log N-O_G(1)$.

\medskip\noindent\textbf{Base case ($d=2$).}
The planar graph bound following Proposition~\ref{prop:spanning} is stronger. Here we prove the bound for $\Phi_2(\mathcal O)$ directly from the energy estimates in Lemma~\ref{lem:d2}, to use as the base of the dimension induction. For $k_i\ge2$, Lemma~\ref{lem:d2}(a) gives the same lower bound $\log N-O_G(1)$ as in dimension one. Suppose $k_i=1$ and put $\tau=\min\{1/6,\varepsilon/(2n)\}$. By Lemma~\ref{lem:d2}(b), all but at most $\tau N$ anchors satisfy $E(a)\le C_\tau N^{4/3+\tau}$. Since $X_i$ is uniform on $P$, histories with such an anchor have total probability at least $1-\tau$. On each of these histories, \eqref{eq:pointwise} and $N-n+1\ge N/2$ for $N\ge2n$ give
\[
H(Y_i\mid\mathcal F_i=h)\ge\log\frac{(N/2)^2}{C_\tau N^{4/3+\tau}}
= (2/3-\tau)\log N-\log(4C_\tau).
\]
On every other history the conditional entropy is nonnegative. Averaging in \eqref{eq:cond}, we obtain
\[
H(Y_i\mid Y_1,\ldots,Y_{i-1})
\ge (1-\tau)(2/3-\tau)\log N-O_{G,\varepsilon}(1)
\ge (2/3-2\tau)\log N-O_{G,\varepsilon}(1).
\]
The steps with $k_i=0$ contribute zero. Summing over all steps and using $2n\tau\le\varepsilon$ gives
\[
H(Y)\ge\bigl(\Phi_2(\mathcal O)-\varepsilon\bigr)\log N-O_{G,\varepsilon}(1).
\]

\medskip\noindent\textbf{Case 1 (non-concentrated, $d\ge3$, $\lambda(P)\le N^\theta$).}
Fix $\eta\in(0,1)$ with $2n\eta\le\varepsilon/2$, and let $i$ be a step with $k_i\ge1$. Apply Theorem~\ref{thm:popstar} with $\eta$ and $\varepsilon/(2n)$. It gives a set $\mathcal B_i\subseteq P^{\underline{k_i}}$ with $|\mathcal B_i|\le\eta N^{k_i}$ such that every $X\in P^{\underline{k_i}}\setminus\mathcal B_i$ satisfies
\[
E(X)\le C\,N^{2-e_{d,k_i}(\theta)+\varepsilon/(2n)},\qquad C=C(G,d,\eta,\varepsilon,\theta).
\]
The root tuple $X_i$ is uniform on $P^{\underline{k_i}}$. This set has size at least $(N-n)^{k_i}\ge N^{k_i}/2$ once $N\ge N_1(G)$. Hence the event $\mathcal E_i=\{X_i\notin\mathcal B_i\}$ satisfies $\mathbb P(\mathcal E_i^c)\le2\eta$. Since $X_i$ is a function of $\mathcal F_i$, the event $\mathcal E_i$ is determined by the history. Using \eqref{eq:cond}, \eqref{eq:pointwise}, and nonnegativity of entropy,
\begin{align*}
H(Y_i\mid Y_1,\ldots,Y_{i-1})
&\ge \sum_{h:\,X_i(h)\notin\mathcal B_i}\mathbb P(\mathcal F_i=h)\,\log\frac{(N-n+1)^2}{E(X_i(h))}\\
&\ge (1-2\eta)\Bigl(\bigl(e_{d,k_i}(\theta)-\tfrac{\varepsilon}{2n}\bigr)\log N-O_G(1)\Bigr).
\end{align*}
Steps with $k_i=0$ contribute $0=e_{d,0}(\theta)$. Summing over $i$ and using $E_d(\mathcal O,\theta)\le n$,
\[
H(Y)\ \ge\ (1-2\eta)\Bigl(E_d(\mathcal O,\theta)-\tfrac{\varepsilon}{2}\Bigr)\log N-O_G(1)
\ \ge\ \Bigl(E_d(\mathcal O,\theta)-\tfrac{\varepsilon}{2}-2n\eta\Bigr)\log N-O_G(1).
\]
Since $2n\eta\le\varepsilon/2$, this gives $|\Delta_G^\circ(P)|\gtrsim N^{E_d(\mathcal O,\theta)-\varepsilon}$. For $N<N_1(G)$, the bound is absorbed into the implicit constant.

\medskip\noindent\textbf{Case 2 (concentrated, $d\ge3$, $\lambda(P)>N^\theta$).}
There exists an affine $(d-2)$-flat $\Gamma$ such that the subset $Q:=P\cap\Gamma$ satisfies $|Q|>N^\theta$. If $\theta\ge\theta_0 = \frac{\log|V(G)|}{\log N}$, then $|Q| > N^\theta \ge |V(G)|$. Thus, there is an injective map $V\to Q$.
Let $\operatorname{aff}(Q)$ be the smallest affine subspace containing $Q$, and write $r=\dim\operatorname{aff}(Q)$ for its dimension. Since $Q$ contains at least two points and lies in the $(d-2)$-flat $\Gamma$, we have $1\le r\le d-2$. Since $r\le d-2<d$, the induction hypothesis (strong induction on the ambient dimension) applies in dimension $r$. Applying it to $Q$ inside 
$\operatorname{aff}(Q)\cong\mathbb R^r$, we obtain:
\[
|\Delta_G^\circ(Q)| \gtrsim |Q|^{\Phi_r(\mathcal O)-\varepsilon}.
\]
By the monotonicity lemma (Lemma~\ref{lem:mono}), $\Phi_r(\mathcal O)\ge\Phi_{d-2}(\mathcal O)$ for any $r\le d-2$. Since any embedding into $Q$ is also an embedding into $P$, $\Delta_G^\circ(Q)\subseteq\Delta_G^\circ(P)$, yielding:
\[
|\Delta_G^\circ(P)| \ge |\Delta_G^\circ(Q)|
\gtrsim |Q|^{\Phi_r(\mathcal O)}|Q|^{-\varepsilon}
\ge N^{\theta\Phi_{d-2}(\mathcal O)-\varepsilon}.
\]

\medskip\noindent\textbf{Conclusion ($d\ge3$).}
Let $\varepsilon>0$ and $\delta>0$. By the definition of $\Phi_d(\mathcal O)$ as a supremum, there is $\theta^*\in[0,1)$, depending only on $G$, $\mathcal O$, $d$ and $\delta$, with
\[
\min\Bigl(E_d(\mathcal O,\theta^*),\ \theta^*\Phi_{d-2}(\mathcal O)\Bigr) \ge \Phi_d(\mathcal O)-\delta.
\]
We may take $\theta^*>0$. If $\delta<\Phi_d(\mathcal O)$, the right-hand side is positive, while the left-hand side is $0$ at $\theta^*=0$. If $\delta\ge\Phi_d(\mathcal O)$, any $\theta^*\in(0,1)$ works.
Because $\theta^* > 0$ is a fixed constant, as $N \to \infty$ the threshold $\theta_0 = \frac{\log |V(G)|}{\log N} \to 0$; hence there exists $N_0$ such that $\theta^* \ge \theta_0$ for all $N \ge N_0$.

For any such $P$ with $N \ge N_0$, we apply the dichotomy at $\theta = \theta^*$. In either case,
\[
|\Delta_G^\circ(P)|\gtrsim_{G,d,\mathcal O,\varepsilon} N^{\min(E_d(\mathcal O,\theta^*),\ \theta^*\Phi_{d-2}(\mathcal O))-\varepsilon} \ge N^{\Phi_d(\mathcal O)-\delta-\varepsilon}.
\]
We have thus shown that for every $\delta>0$ and every $\varepsilon>0$,
\[
|\Delta_G^\circ(P)|\gtrsim_{G,d,\mathcal O,\delta,\varepsilon} N^{\Phi_d(\mathcal O)-\delta-\varepsilon}.
\]
Given any $\varepsilon'>0$, applying this with $\delta=\varepsilon=\varepsilon'/2$ yields $|\Delta_G^\circ(P)|\gtrsim N^{\Phi_d(\mathcal O)-\varepsilon'}$, as required. This completes the proof.
\end{proof}

\section{Main theorem}\label{sec:main}

Throughout this section, $\mathcal O$ is a $k$-degenerate construction order with $m=\sum_i k_i\ge1$, and we write $\Phi_r=\Phi_r(\mathcal O)$. By the proof of Lemma~\ref{lem:deg}, $\Phi_r>0$ for every $r\ge1$, so we may set
\[
u_r:=\frac{m}{\Phi_r}\qquad(r\ge1).
\]
In this notation, Lemma~\ref{lem:deg} states that
\begin{equation}\label{eq:u-recursion}
u_r=\max\bigl(r-1,\ u_{r-2}+1\bigr)\qquad\text{for every } r\ge\max(k+1,3).
\end{equation}

\begin{lemma}\label{lem:iterate}
Let $b\ge1$ be an integer with $b+2\ge\max(k+1,3)$. Then for every $j\ge1$,
\[
u_{b+2j}=\max\bigl(b+2j-1,\ u_b+j\bigr).
\]
In particular, for $j\ge1$ we have $\Phi_{b+2j}=\frac{m}{b+2j-1}$ if and only if $j\ge u_b-b+1$.
\end{lemma}
\begin{proof}
For $j\ge1$, the index $r=b+2j$ satisfies $r\ge b+2\ge\max(k+1,3)$, so \eqref{eq:u-recursion} applies at $r$. We induct on $j$. For $j=1$, \eqref{eq:u-recursion} gives $u_{b+2}=\max(b+1,\ u_b+1)$. If the formula holds for some $j\ge1$, then
\[
u_{b+2j+2}=\max\bigl(b+2j+1,\ u_{b+2j}+1\bigr)=\max\bigl(b+2j+1,\ b+2j,\ u_b+j+1\bigr)=\max\bigl(b+2(j+1)-1,\ u_b+(j+1)\bigr).
\]
Finally, $u_{b+2j}=b+2j-1$ holds exactly when $u_b+j\le b+2j-1$, that is, when $j\ge u_b-b+1$.
\end{proof}

\begin{lemma}\label{lem:u-bound}
If $k\ge2$, then $u_r\le k+\lfloor\frac{r-1}{2}\rfloor$ for every $1\le r\le k+2$.
\end{lemma}
\begin{proof}
We argue by induction on $r$, starting with $r=1,2$ and using the bound at $r-2$ in the inductive step. Let $n_j=|\{i:k_i=j\}|$.

\medskip\noindent\textbf{Base case ($r=1,2$).}
Since $m=\sum_{j=1}^k jn_j\le k\sum_{j=1}^k n_j=k\Phi_1$, we have $u_1\le k$.
For $r=2$, $\Phi_2=\frac23n_1+\sum_{j=2}^k n_j$. We compare this term by term with $\frac mk=\frac{n_1}{k}+\sum_{j=2}^k\frac jk n_j$. Since $k\ge2$, we have $\frac23\ge\frac12\ge\frac1k$, and $1\ge\frac jk$ for $j\le k$. Thus $\Phi_2\ge\frac mk$, that is, $u_2\le k$.

\medskip\noindent\textbf{Inductive step ($3\le r\le k+2$).}
Fix $3\le r\le k+2$ and $\theta\in[0,1)$, and define
\[
f_r(x)=\min\Bigl(1,\,x(1-\theta),\,\frac{x}{r-1}\Bigr),\qquad x\ge0.
\]
For fixed $\theta$, $f_r$ is the minimum of three linear functions of $x$, so it is concave; since also $f_r(0)=0$, concavity applied to $x=\frac xk\cdot k+\bigl(1-\frac xk\bigr)\cdot0$ gives $f_r(x)\ge\frac{x}{k}f_r(k)$ for all $x\in[0,k]$. Moreover, $e_{r,x}(\theta)\ge f_r(x)$ for every integer $x\ge0$: for $x\le r-1$ the two coincide, and for $x\ge r$ the convention gives $e_{r,x}(\theta)=1\ge f_r(x)$. Applying this to each $k_i\le k$,
\[
E_r(\theta)=\sum_{i=1}^n e_{r,k_i}(\theta)\ \ge\ \sum_{i=1}^n\frac{k_i}{k}\,f_r(k)=m\min\Bigl(\frac1k,\,1-\theta,\,\frac1{r-1}\Bigr).
\]
Write $c=\min\bigl(\frac mk,\frac m{r-1}\bigr)$ and $s=\frac{m\Phi_{r-2}}{m+\Phi_{r-2}}$. Then
\[
\Phi_r\ \ge\ \sup_{0\le\theta<1}\min\bigl(c,\ m(1-\theta),\ \theta\Phi_{r-2}\bigr)\ =\ \min(c,s),
\]
since $\min\bigl(m(1-\theta),\theta\Phi_{r-2}\bigr)\le s$ for every $\theta$, with equality at $\theta=\frac{m}{m+\Phi_{r-2}}\in(0,1)$. Taking reciprocals,
\[
u_r\le\max\bigl(k,\ r-1,\ u_{r-2}+1\bigr).
\]
By the induction hypothesis and $\lfloor\frac{r-3}{2}\rfloor+1=\lfloor\frac{r-1}{2}\rfloor$, this gives
\[
u_r\le\max\Bigl(k,\ r-1,\ k+\Bigl\lfloor\frac{r-1}{2}\Bigr\rfloor\Bigr).
\]
Clearly $k\le k+\lfloor\frac{r-1}{2}\rfloor$. Also, $r-1\le k+\lfloor\frac{r-1}{2}\rfloor$ is equivalent to $\lceil\frac{r-1}{2}\rceil\le k$, which holds because $\lceil\frac{r-1}{2}\rceil\le\lceil\frac{k+1}{2}\rceil\le k$. Hence $u_r\le k+\lfloor\frac{r-1}{2}\rfloor$.
\end{proof}

\begin{proposition}\label{prop:threshold}
$\Phi_d(\mathcal O)=\frac{m}{d-1}$ for every $d\ge\max(2k,3)$.
\end{proposition}
\begin{proof}
\medskip\noindent\textbf{Case $k=1$.}
Every step has $k_i\in\{0,1\}$, so $\Phi_1=m$ and $\Phi_2=\frac23m$, that is, $u_1=1$ and $u_2=\frac32$. Lemma~\ref{lem:iterate} with $b=1$ applies, since $b+2=3=\max(k+1,3)$; as $u_1-1+1=1$, it gives $\Phi_d=\frac m{d-1}$ for every odd $d\ge3$. Lemma~\ref{lem:iterate} with $b=2$ also applies, since $b+2=4\ge3$; as $u_2-2+1=\frac12\le1$, it gives the same for every even $d\ge4$.

\medskip\noindent\textbf{Case $k\ge2$.}
We apply Lemma~\ref{lem:iterate} with $b=k$ and with $b=k-1$. Its hypothesis holds in both cases: for $b=k$ we have $b+2=k+2\ge\max(k+1,3)$, and for $b=k-1$ we have $b\ge1$ and $b+2=k+1=\max(k+1,3)$, since $k\ge2$.

For $d\equiv k\pmod 2$: by Lemma~\ref{lem:u-bound} with $r=k$, we have $u_k-k+1\le\lfloor\frac{k-1}{2}\rfloor+1$. Hence Lemma~\ref{lem:iterate} with $b=k$ gives $\Phi_d=\frac m{d-1}$ for every such $d$ with
\[
d\ \ge\ k+2\Bigl(\Bigl\lfloor\frac{k-1}{2}\Bigr\rfloor+1\Bigr)=\begin{cases}2k+1,&k\text{ odd},\\ 2k,&k\text{ even}.\end{cases}
\]

For $d\not\equiv k\pmod 2$: by Lemma~\ref{lem:u-bound} with $r=k-1$, we have $u_{k-1}-(k-1)+1\le\lfloor\frac{k-2}{2}\rfloor+2=\lfloor\frac k2\rfloor+1$. Hence Lemma~\ref{lem:iterate} with $b=k-1$ gives $\Phi_d=\frac m{d-1}$ for every such $d$ with
\[
d\ \ge\ (k-1)+2\Bigl(\Bigl\lfloor\frac k2\Bigr\rfloor+1\Bigr)=\begin{cases}2k,&k\text{ odd},\\ 2k+1,&k\text{ even}.\end{cases}
\]
In both cases, every even $d\ge2k$ and every odd $d\ge2k+1$ is covered, that is, every $d\ge2k$.
\end{proof}

\begin{theorem}\label{thm:main}
Let $G=(V,E)$ be a finite graph with $m\ge1$ edges, and suppose $G$ is $k$-degenerate for some $k\ge1$. Then for every $d\ge\max(2k,3)$, every finite $P\subset\mathbb R^d$ with $|P|=N\ge|V|$, and every $\varepsilon>0$,
\[
|\Delta_G^\circ(P)|\gtrsim_{G,d,\varepsilon} N^{\frac{m}{d-1}-\varepsilon}.
\]
\end{theorem}
\begin{proof}
Fix a $k$-degenerate construction order $\mathcal O$ of $G$. By Proposition~\ref{prop:threshold}, $\Phi_d(\mathcal O)=m/(d-1)$ for every $d\ge\max(2k,3)$, and Theorem~\ref{thm:entropy} gives $|\Delta_G^\circ(P)|\gtrsim N^{m/(d-1)-\varepsilon}$.
\end{proof}

\section*{Acknowledgments}
The authors are grateful to Dora Woodruff for helpful conversations on the project. The authors thank the organizers of the 2026 Polymath Jr.\ program, during which this research was carried out. The program is supported in part by the National Science Foundation under award DMS-2341670.

\section*{AI Statement}
The authors acknowledge the use of GPT-5.5 Pro and GPT-5.6 Sol Pro (OpenAI) during the initial exploration of the problem. These models assisted in particular in the discovery of Theorem~\ref{thm:popstar}, and provided limited assistance in the discovery of other results. All statements and proofs were subsequently verified, substantially revised, and reformulated by the authors, who take full responsibility for the content of this paper.

\bibliographystyle{plain}
\bibliography{references}

\end{document}